\documentclass[9pt,technote]{IEEEtran}
\usepackage{ifpdf}
\usepackage{cite}

\ifCLASSINFOpdf
 \usepackage[pdftex]{graphicx}
\else
 \usepackage[dvips]{graphicx}
\fi
\usepackage[cmex10]{amsmath}
\usepackage{amsfonts}
\usepackage{array}
\usepackage{fixltx2e}
\usepackage{url}
\usepackage{textcomp}
\usepackage{algorithm}
\usepackage{algorithmic}

\newtheorem{thm}{Theorem}[section] \newtheorem{cor}[thm]{Corollary}
\newtheorem{lem}[thm]{Lemma} 
 
\newtheorem{ex}[thm]{Example} \newtheorem{defn}[thm]{Definition}

\newcommand{\R}{\mathbb R} \newcommand{\N}{\mathbb N}
\newcommand{\n}{\underline{n}}
\newcommand{\m}{\underline{m}}
\newcommand{\eps}{\varepsilon}
\newcommand{\norm}[1]{\left\Vert#1\right\Vert}
\newcommand{\abs}[1]{\left\vert#1\right\vert}

\DeclareMathOperator*{\conv}{conv}
\newcommand{\yconv}{\conv_{i\in\m}y^i}

\begin{document}
\title{On Conditions for Convergence to Consensus}

\author{Jan Lorenz, Dirk A. Lorenz% <-this % stops a space
  \thanks{J. Lorenz is with the Chair of Systems Design, ETH Zurich,
    Kreuzplatz 5, 8032 Zurich, Switzerland, most of the work was done
    when he was with the Department of Mathematics and Computer
    Science, University of Bremen, Bibliothekstrasse 1, 28359 Bremen,
    Germany, post@janlo.de,
    {http://www.janlo.de.}}% <-this % stops a space
  \thanks{D.A. Lorenz is with the Institute for Analysis and Algebra,
    Carl-Friedrich Gau\ss\ Department, TU Braunschweig,
    38092~Braunschweig, Germany, d.lorenz@tu-braunschweig.de,
    http://www.tu-braunschweig.de/iaa/personal/lorenz.}% <-this % stops a space
  \thanks{Manuscript received xxx 00, 2008; revised xxx 00, 2008.}}

\markboth{}%
{ }
% The only time the second header will appear is for the odd numbered
% pages after the title page when using the twoside option.
% *** name in the headers of peer review papers.  *** You can use
% \ifCLASSOPTIONpeerreview for conditional compilation here if you
% desire.

% this: \IEEEpubid{0000--0000/00\$00.00~\copyright~2007 IEEE} Remember,
% if you use this you must call \IEEEpubidadjcol in the second column for
% its text to clear the IEEEpubid mark.

% use for special paper notices \IEEEspecialpapernotice{(Invited Paper)}
\maketitle

\begin{abstract}
  % \boldmath
  A new theorem on conditions for convergence to consensus of a
  multiagent time-dependent time-discrete dynamical system is presented.
  The theorem  is build up on the notion of averaging maps. We compare this theorem to results by
  Moreau (IEEE Transactions
  on Automatic Control, vol. 50, no. 2, 2005) about set-valued Lyapunov theory and convergence under switching communication topologies. We give examples
  that point out differences of approaches including examples where
  Moreau's theorem is not applicable but ours is. Further on, we give
  examples that demonstrate that the theory of convergence to
  consensus is still not complete.
\end{abstract}


% For peer review papers, you can put extra information on the cover page
% as needed: \ifCLASSOPTIONpeerreview
% \begin{center} \bfseries EDICS Category: 3-BBND \end{center} \fi
\IEEEpeerreviewmaketitle

\section{Introduction}
\label{sec:introduction}
In this technical note we analyze discrete dynamical systems of consensus
formation as presented in the context of distributed computing
\cite{Tsitsiklis1984ProblemsinDecentralized,Tsitsiklis.Bertsekas.ea1986Distributedasynchronousdeterministic}, flocking (e.g. of unmanned
aerial vehicles)
\cite{Jadbabaie.Lin.ea2003Coordinationofgroups,Saber.Murray2003Flockingwithobstacle,Blondel.Hendrickx.ea2005ConvergenceinMultiagent} and
general as multi-agent coordination problems
\cite{Moreau2005Stabilityofmultiagent,Ren.Beard.ea2005surveyofconsensus,Xiao.Boyd2004Fastlineariterations} (to mention just a
few). The dynamical system may also be called `agreement algorithm' or
`consensus protocol'. The convergence theorems of Moreau \cite{Moreau2005Stabilityofmultiagent}
together with the extensions of Angeli and Bliman
\cite{Angeli.Bliman2006Stabilityofleaderless} are the most general ones.
The main theorem of Moreau states conditions for
convergence to consensus under switching communication topologies.
Convergence to consensus is there implied by `global asymptotic stability of
the set of equilibrium solutions with consensus as equilibrium points'.
%($x\in S^n$ is a consensus if $x^1 = \dots = x^n$).
Conditions are on
the one hand on the communication topologies in their time-evolution and
on the other hand on the updating maps. Moreau applied a set-valued
Lyapunov theory, which uses a set-valued function on the state space
which is contractive with respect to the updating map. This implies
convergence of the set to a singleton.

We contribute a similar but new approach based on the notion of an
averaging map. Moreau deals with communication topologies by defining
conditions on how many successive communication topologies must be
regarded until the composition of these updating maps fulfills the
contraction properties used to apply the set-valued Lyapunov theory. We
skip the issue on changing communication topologies and deal directly
with maps which fulfill a contraction property which is different from
Moreau's.

Our theorem generalizes a result of Krause \cite{Krause2009Compromiseconsensusand} by allowing
arbitrary switching between different averaging maps but follow the same
line of compactness, continuity and convexity arguments.  

Section \ref{sec:convergence-result} presents the convergence result
and possible extensions. Section \ref{sec:comp-set-valu} discusses the
relations to two of Moreau's theorems in more
detail. Section \ref{sec:counterexamples} gives examples and
counterexamples to show existing gaps in the theory of consensus
algorithms.  All proofs of lemmas and theorems are collected in
Appendix~\ref{app:proof-theor-familyavmap}.

\section{Convergence result}
\label{sec:convergence-result}
We consider a dynamical system of the form
\begin{equation}
x(t+1) = f_t(x(t))\label{eq:dynsys}
\end{equation}
with discrete time $t\in\N$. Dynamics take place
in a $d\times n$-dimensional space:  We consider a set of
agents $\n = \{1,\dots,n\}$ where each of them has coordinates in a $d$-dimensional set $S\subset \R^d$.  Hence, the solutions of (\ref{eq:dynsys}) have the form $x:\N\to S^n\subset \R^{d\times n}$.
The individual
coordinates of agent $i$ at time $t\in\N$ is labeled $x^i(t) \in S$, and
$x(t)\in S^n$ is called the \emph{profile} at time $t\in\N$.
Finally, the mappings $f_t$ which govern the dynamics are of the
form $f_t:S^n\to S^n$. We denote the component functions by $f_t^i$.

To state our main result on convergence of such systems to consensus we
introduce the following notations.  An element $x\in S^n$ is called
\emph{consensus} if all $d$-dimensional coordinates $x^i$ have the same
value, i.e.~there exists a vector $\gamma\in S$ such that $x^i = \gamma$
for $i\in \n$.  By $\conv_{i\in\n} x^i$ we define the convex hull of
the vectors $x^1,\dots,x^n$.

The core notion in this note is an `averaging map'. We build the
definition of an averaging map on a generalized convex hull. Consider a
continuous function $y: S^n \to S^m$ which maps a profile to a
certain set of $m$ vectors $y(x)=(y^1(x),\dots,y^m(x))$ such that for all
$x\in S^n$ and all $i\in\n$ it holds $x^i \in \conv_{j\in\m}y^j(x)$. We
call such a function $y$ a generalized \emph{barycentric coordinate map}
and we call $\conv_{j\in\m}y^j(x)$ the \emph{$y$-convex hull}
of the vectors $x^1,\dots,x^n$. (We call $y$ `generalized' because it needs not be
a bijective transformation.) So, a $y$-convex hull is a set-valued
function from $S^n$ to the compact and convex subsets of $S$. 
We call a set $y$-convex, if it is the union of the $y$-convex hulls of 
all $n$ of its points. Examples
for $y$-convex hulls include the convex hull itself, and the
multidimensional interval $[\min_{i\in\n}x^i,\max_{i\in\n}x^i]$ (with
$\min$ and $\max$ applied componentwise).  For the first it
holds $m=n$ for the second $m = 2^d$. Many other examples fit into this
setting: the smallest interval for any basis of $\R^d$ \cite[Example
2]{Angeli.Bliman2006Stabilityofleaderless}, or smallest polytope with
faces parallel to a set of $k\geq d+1$ hyperplanes \cite[Example
3]{Angeli.Bliman2006Stabilityofleaderless} containing $x^1,\dots,x^n$ (the generalized barycentric coordinates are then
the extreme points of the polytope, perhaps with multiples to have a
constant $m$).  %\todo{Bild zu $y$-avmap machen}
Now, we define the central notion of this paper.
\begin{defn}
  \label{defn:avmap}
  Let $S\subset \R^d$, $y:S^n\to S^m$ be a generalized barycentric
  coordinate map such that $S$ is $y$-convex. A mapping $f:S^n\to S^n$ is called a
  \emph{$y$-averaging map}, if for every $x\in S^n$ it holds
  \begin{equation}
    \label{eq:avconv}
    \yconv(f(x))\subset\yconv(x).
  \end{equation}
  Furthermore, a \emph{proper $y$-averaging map} is a $y$-averaging map, such
  that for every $x\in S^n$ which is not a consensus, the above inclusion
  is strict.
\end{defn}
 
A $y$-averaging map maps a profile $x$ into its $y$-convex hull.
Furthermore, the $y$-convex hull of the new
profile $f(x)$ lies in the $y$-convex hull of the vectors $x^1,\dots,x^n$.
Hence, we may also work with the $y$-convex hull of the initial profile $x(0)$ instead of the set $S$.
%(Notice that this
%assumption is automatically fulfilled in all the examples mentioned in
%the former paragraph if we claim for all $i\in\n$ that $x^i \in
%\yconv (x)$.)
Sometimes it is useful to look at the
contraposition of the definition of proper: If equality holds in
\eqref{eq:avconv} this implies that $x$ is a consensus. In the following we may omit
'$y$' when we mention an averaging map, but for an averaging map the
definition of $y$ is a prerequisite. The best proxy for the mind is $y = \textrm{id}$.

Since we are going to consider families of averaging maps we introduce
the concept of equiproper averaging maps. To this end, we need the
Hausdorff distance on the set of compact subsets of a metric space $(X,d)$.
The distance of a point $x\in X$ and a nonempty compact set $C\subset X$
is defined as $d(x,C) := \min_{c \in C} d(x,c)$.
Let $B,C \subset X$ be nonempty and compact, then the \emph{Hausdorff
  distance} is defined as
\[
d_H(B,C) := \max\{ \max_{b \in B}d(b,C) , \max_{c \in C}d(c,B) \}.
\]
Equivalently, one can say that the Hausdorff distance is the smallest $\eps$
such that the $\eps$-neighborhood of $B$ contains $C$ and the
$\eps$-neighborhood of $C$ contains $B$. It is easy to see that $d_H(B,C)
= 0$ holds if and only if $B=C$.  In the special case $B\subset C\subset S\subset \R^d$ it holds 
\begin{equation}
  \label{eq:d_H_inclusion}
  d_H(B,C) = \max_{b \in B}d(b,C) = \max_{b\in B}\min_{c\in C} \norm{b-c}.
\end{equation}

\begin{defn}
  \label{defn:equiproper}
  Let $y$ be a generalized barycentric coordinate map and let $F$
  be a family of proper $y$-averaging maps.  $F$ is called
  \emph{equiproper}, if for every $x\in S^n$ which is not a consensus,
  there is $\delta(x)>0$ such that for all $f\in F$
  \begin{equation}
    d_H\bigl(\yconv(f(x)),\yconv(x)\bigr) > \delta(x).
    \label{eq:convequiproper}
  \end{equation}
\end{defn}
Now we state a lemma which says that the family of equiproper
$y$-averaging maps is closed under pointwise limits.
\begin{lem}
  \label{lem:limit_proper_avmap}
  Let $f_t$ be a sequence of $y$-averaging maps forming an equiproper
  family of $y$-averaging maps such that $f_t\to g$ pointwise. Then $g$
  is a proper $y$-averaging map.
\end{lem}
Now we are able to state our main theorem.
\begin{thm}\label{thm:familyavmap}
  Let $S\subset \R^d$, $y$ be a
  generalized barycentric coordinate map such that $S$ is $y$-convex,
  and $F$ be an
  equicontinuous family of equiproper $y$-averaging maps on $S^n$. Then it
  holds for any sequence $(f_t)_{t\in\N}$ with $f_t\in F$ and any $x(0)
  \in S^n$ that
  the solution of~(\ref{eq:dynsys}) converges to a consensus,
  i.e. there exists $\gamma\in S$ such that for all $i\in \n$ it holds
  $\lim_{t\to\infty} x^i(t) = \gamma$.
\end{thm}

Notice that the limit $\gamma$ depends not only on the initial value
$x(0)$ but also on the realization of the sequence $(f_t)_{t\in\N}$,
however, $\gamma$ depends continuously on the intial value if the
sequence $(f_t)$ is fixed as the following lemma and corollary show.
\begin{lem}\label{lem:cont_dep_initial_value}
  Let $(X,d)$ be a metric space and $f_t:X\to X$ be such that the solution
  of $x(t+1) = f_t(x(t))$ converge to some limit for every initial
  value $x(0)\in X$. Then the limit depends continuously on the initial
  value if $\{f_t\}$ is an equicontinuous family.
\end{lem}
The following corollary is a direct consequence.
\begin{cor}
  Let the sequence $(f_t)$ in the situation of
  Theorem~\ref{thm:familyavmap} be fixed. Then the
  consensus value $\gamma$ (which exists due to Theorem~\ref{thm:familyavmap}) depends continuously 
  on the initial value.
\end{cor}

Theorem~\ref{thm:familyavmap} is a generalization of a theorem of Krause~\cite{Krause2009Compromiseconsensusand}.
Krause's theorem is the special case when $y$ is the identity and $F$
contains only one proper averaging map. Notice that 'equi' in equiproper
and equicontinuous can be omitted if $F$ is a finite set. An easy
extension is to allow $F$ to contain also non-proper averaging maps (but
at least one proper averaging map). Then the sequence $(f_t)_{t\in\N}$ has to contain
a subsequence $(f_{t_s})_{s\in\N}$ of equiproper averaging maps to ensure
convergence to
consensus. This holds because then $\{g_s \,|\, g_s = f_{t_s}\circ\dots\circ f_{t_{s+1}}\}$
is an equiproper set of averaging maps for $s\in\N$. Notice that it is possible that a sequence of averaging
maps contains a subsequence as above such that subcompositions $g_s$ form
an equiproper set, even when no $f_t$ is proper. The easiest example is when $F$ contains
only one linear map which is determined by a row-stochastic square matrix
which is regular but not scrambling (see Seneta \cite{Seneta2006Non-NegativeMatricesand}). For
linear systems `row-stochastic' is equivalent to `being an averaging map' (with $y$ the
identity) and `scrambling' is equivalent to `proper'. From the theory of
nonnegative matrices we know that for each regular matrix there is an
integer such that higher powers are scrambling.

In the spirit of~\cite{Angeli.Bliman2006Stabilityofleaderless} we
state another generalization of Theorem~\ref{thm:familyavmap} which
deals with deformations of the hull.  To this end, let $S,T\subset
\R^d$ be compact and $\phi:T\to S$ be a homeomorphism. For a
generalized barycentric coordinate map $y:S^n\to S^m$ we define the
$y,\phi$-hull as $\phi^{-1}( \yconv(\phi(x)))$.  Now, a
$y,\phi$-averaging map $g$ is defined analogous to
Definition~\ref{defn:avmap}:
\[
\phi^{-1}(\conv_{i\in\m}y^i(\phi(g(x))))\subset \phi^{-1}(\conv_{i\in\m} y^i(\phi(x))).
\]
Note, that the $y,\phi$-hull is not necessarily convex,
see\cite[Example 6]{Angeli.Bliman2006Stabilityofleaderless}.
The extension of the notions `proper' and `equiproper' is straightforward. 

\begin{thm}
  \label{thm:family_yconvphi_avmap}
  Let $\phi:T\to S$ be continuous with Lipschitz continuous inverse
  and let $y$ be a generalized barycentric coordinate map such that $S$ is $y$-convex.
  Let $G$ be a family of equicontinuous, equiproper $y,\phi$-averaging
  maps on $T^n$.  Then it holds for any sequence $(g_t)_{t\in\N}$ with $g_t\in G$
  and any $x(0) \in T^n$ that the solution of $x(t+1) = g_t(x(t))$
  converges to a consensus.
\end{thm}

\section{Comparison with Moreau's set-valued Lyapunov theory and main
theorem}
\label{sec:comp-set-valu}

Theorem \ref{thm:familyavmap} has similarities to Moreau's set-valued
Lyapunov Theorem~\cite[Theorem 4]{Moreau2005Stabilityofmultiagent}. This theorem implies global asymptotic stability of the set of
equilibrium solutions when there exists a set-valued function $V$ on the
state space, a measure for these sets $\mu$, and a positive definite
function $\beta$ on the state space. Essentially it has to hold
$V(f_t(x))\subset V(x)$ and $\mu(V(f_t(x)))-\mu(V(x)) \leq -\beta(x)$.
%(This is just the essence of the theorem omitting some details.)
The best example to imagine is $V = \conv$, and $\mu$ is the diameter of a set.

The set of equilibrium solutions for the dynamical system \eqref{eq:dynsys} under the conditions of Theorem \ref{thm:familyavmap} contains only all constant solutions on consensus vectors, due to the equiproperness of $F$. Given this set of equilibrium solutions, ``global asymptotic stability of the set of equilibrium solutions'' implies convergence to consensus. Convergence to consensus is thus a special case of the set-valued Lyapunov Theorem in \cite{Moreau2005Stabilityofmultiagent}. To the best of our knowledge, it is the only case in which the theorem has been used so far.

Compared with our Theorem \ref{thm:familyavmap} the role of the set-valued map $V$ is taken by the
$y$-convex hull. So, we also deal with a general class of functions due to the
various possible coordinate maps $y:S^n\to S^m$---we only assume that $m$ is finite.
However, we do not need a general
measure $\mu$ on these maps. The assumption $\mu(V(f_t(x)))-\mu(V(x))
< \beta(x)$ corresponds to $d_H(V(f_t(x)),V(x)) > \delta(x)$. This is
a different condition and often weaker, as for example in the case where Moreau specifies it to proof his main Theorem \cite[Theorem 2]{Moreau2005Stabilityofmultiagent}. There $\mu$ is the diameter of $V(x)$ (which he specifies as the $\conv(x)$).

Theorem \ref{thm:familyavmap} has also similarities to Moreau's main
theorem \cite[Theorem 2]{Moreau2005Stabilityofmultiagent}. This theorem is more specific than Theorem~\ref{thm:familyavmap} by incorporating
switching communication topologies. Its main drawback is that it relies
very much on convex hulls (see \cite{Angeli.Bliman2006Stabilityofleaderless} for a method to
overcome this drawback). Our result generalizes to convex hulls of
generalized coordinate maps.  Further on, in Moreau's theorem agents are
forced to move into the relative interior of the convex hull (respecting
the communication topology). Specifically, this implies that agents have
to leave all extreme points of the convex hull (of agents in its neighborhood) after one iteration. Our
theorem needs only agents at one arbitrary extreme point (of the global $y$-convex
hull) to leave it towards the interior after one iteration. This is implied by properness of averaging maps. 
The assumption `equiproper' in our theorem finds its analog in Moreau's
theorem by assuming that the sets $e_k(\mathcal{A}(t))(x)$ are chosen independently of $t$.\footnote{Here the matrix $\mathcal{A}(t)$ is the arbitrarily chosen communication topology at time $t$ and $x$ is a given state. The set $e_k(\mathcal{A}(t))(x)$ is a subset of the relative interior of the convex hull of the neighbors of $k$ (including $k$) in the current communication topology, and it determines the set where the state of node $k$ has to remain in after one iteration. So, $e_k$ has to be fixed for a given communication topology and a certain state regardless of the chosen updating map $f(t,\cdots)$. This is in analogy to equiproper which implies the existence of a minimal Hausdorff distance $\delta(x)$ after one iteration for a given state but all possible averaging maps.}

Summarizing the above one can say that both Moreau's theorem and
Theorem~\ref{thm:familyavmap} are similar. However, the assumptions as
well as the methods of proof are different. On the one hand we do not
incorporate switching communication topologies explicitly, but on the
other hand we need weaker conditions for the updating maps $f_t$. Further on, we generalized to
$y$-convex hulls and are also able to incorporate the extensions of
Moreau's theorem by Angeli and Bliman
\cite{Angeli.Bliman2006Stabilityofleaderless} to overcome the restriction to convex sets. Moreover, the notion of
a (equi-)proper $y$-averaging map allows a systematic and structured
treatment of consensus algorithms (see e.g.~the results in
Lemma~\ref{lem:limit_proper_avmap} and
Lemma~\ref{lem:cont_dep_initial_value}). 
Hence,
Theorem~\ref{thm:familyavmap} together with
\ref{thm:family_yconvphi_avmap} give an alternative approach to the
analysis of consensus protocols whose applicability is illustrated by
examples in the next section.

\section{Examples and Counterexamples}
\label{sec:counterexamples}
In this section we present counterexamples (Examples \ref{ex:avmap3}--\ref{ex:vanishing}) to point that the existing theory, including our Theorem~\ref{thm:familyavmap}, delivers no sharp results on convergence to consensus. We also give examples which show cases, where our theorem is applicable but Theorem 2 of Moreau is not (Examples \ref{ex:watergun}--\ref{ex:nonarithmetic}).

Continuity, for instance, is not necessary for convergence to consensus
since there are discontinuous proper averaging maps which converge to consensus
(one may take different averaging maps on different subdomains of $S$).
On the other hand discontinuity may destroy convergence to consensus even
for proper averaging maps (see\cite[Section 3.1]{Lorenz2007RepeatedAveragingand} for examples for this phenomenon).

The next two examples illustrate the role of equiproperness.
\begin{ex}[Non-equiproper not leading to consensus]\label{ex:avmap3}
  Let %$f_t: (\R)^2 \to (\R)^2$ with
  \[
  f_t(x^1,x^2) := \left((1-\frac{1}{4^t})x^1+\frac{1}{4^t}x^2 \ , \
    \frac{1}{4^t}x^1+(1-\frac{1}{4^t})x^2\right)
  \]
  It is easy to see that for $t\geq 1$ and $x(1) = (0,1)$ it holds that
  $x^1(t) < \frac{1}{3}$ and $x^2(t) > \frac{2}{3}$.
  Obviously, $\{f_t \,|\, t\in\N\}$ is not equiproper because $f_t$
  converges to the identity as $t\to\infty$.
\end{ex}

\begin{ex}[Non-equiproper leading to consensus]\label{ex:notequiconv}
  Let %$f_t: (\R)^2 \to (\R)^2$ with
  \[
  f_t(x^1,x^2) := \left((1-\frac{1}{t})x^1+\frac{1}{t}x^2 \ , \
    x^2\right)
  \]
  This example is not equiproper, because $f_t$ converges to the identity
  for $t\to\infty$. Thus, Theorem~\ref{thm:familyavmap} does not apply, but for $t\geq 2$ and any $x(2) \in (\R)^2$ the
  system $x(t+1)=f_t(x(t))$ has the solution $x(t) = (\frac{1}{t-1}x^1(2)
  + \frac{t-2}{t-1}x^2(2)\ ,\, x^2(2))$ and thus converges to consensus at
  $x^2(2)$.
  Note that the convergence is not at an exponential rate.
\end{ex}
Convergence to consensus in the last example can also not be ensured by
Moreau's theorems. 

The next example illustrates the role of equicontinuity and is
inspired by bounded confidence~\cite{Hegselmann.Krause2002OpinionDynamicsand}.

\begin{ex}[Vanishing confidence]\label{ex:vanishing}
  Let $f_t: \R^n \to \R^n$ with
  \[
  f_t^i(x) := \frac{\sum_{j=1}^n D_t(|x^i - x^j|)x^j}{\sum_{j=1}^n
    D_t(|x^i - x^j|)}
  \]
  and $D_t: \R_{\geq 0} \to \R_{\geq 0}$. Now, $f_t$ is an averaging map
  for any choice of $D_t$. Further on, $f_t$ is continuous if $D_t$ is,
  and $f_t$ is proper if $D_t$ is strictly positive.  We chose $ D_t(y)
  := e^{-(\frac{y}{\eps})^t} $ as a sequence of functions which has the
  cutoff function as pointwise limit function. Hence, $D_t$ is continuous
  but $\{D_t\,|\,t\in\N\}$ is not equicontinuous.  For $x(0) = (0,8),$
  $\eps=1$ the process $x(t)=f_t(x(t))$ does not converge to consensus
  although only proper averaging maps are involved.  Rough estimates show
  that $\abs{x^1(t) - x^2(t)} \geq 4$.

  For other settings convergence under vanishing confidence is possible,
  as numerical examples in \cite{Lorenz2007RepeatedAveragingand} show.
\end{ex}

The following examples are to show limitations of Moreau's Theorem 2 and how Theorem~\ref{thm:familyavmap} can be applied to show convergence to consensus.
\begin{ex}[Rendezvous problem with watergun sensors] \label{ex:watergun}
  We consider a version of the Rendezvous Problem
  \cite{Lin.Morse.ea2003multi-agentrendezvousproblem} where $n$
  agents are to locate themselves decentralized at the
  same position in twodimensional space. Each agent has three waterguns, an \emph{activation gun} and two \emph{search guns}. Agents can perceive from which kind of gun they were hit and can respond (e.g. acoustically). The search gun is used as a sensor to check if there is at least one other agent in direction $\alpha\in[0,2\pi[$. The activation gun is used to activate other agents. When another agent responds to a shot by the activation gun, the shooting agent switches to standby (only responding if hit). With two search guns an agent can particularily perform a move into direction $\beta \in [0,2\pi[$ under \\
\emph{Rule ($\ast$):} Move until either the position of an other agent is reached or until there is an agent in the directions $\beta+\tfrac{\pi}{2}$ or $\beta-\tfrac{\pi}{2}$. (Move while constantly shooting left and right with search gun until someone is hit.) \\
Initially the $n$ agents are located at different positions in space and the multi-agent protocol is started form the outside by activating one agent. Whenever an agent is activated it executes the following program:
\begin{algorithmic}
 \STATE search gun all around shot, detect $A$ as set of all directions where agents are
 \STATE select $\alpha,\gamma$ such that for all $-1\leq c\leq1$ it holds $(\alpha+c\gamma) \mod 2\pi \notin A$ and $\gamma$ maximal
 \IF {$\gamma \geq \frac{\pi}{2}+\frac{\pi}{n}$}
	\STATE tie agents at same position to move together
	\STATE move direction $\beta = \alpha+2\pi \mod 2\pi$ with rule ($\ast$)
 \ENDIF
 \STATE activation gun all around shot (random start) until someone hit
 \IF {no one hit}
	\STATE give signal `consensus found!'
 \ENDIF
\end{algorithmic}
\end{ex}
The protocol ensures that always only one agent is activated when an agent finishes its action unless consensus is found.
It also always leads to the movement of an agent after some time unless consensus is found, because for every configuration there is always at least one agent whose position is an extreme point of the convex hull such that the exterior angle of the convex hull is larger than $\pi+\frac{2\pi}{n}$ and thus $\gamma\geq\frac{\pi}{2}+\frac{\pi}{n}$. This is because $\pi+\frac{2\pi}{n}$ is the exterior angle of a regular polygon with $n$ edges, which is the `worst case'-polygon. It is `worst case', because it has from all polygons with $n$ edges the largest minimal exterior angle. Thus, the random search for an agent which finds a direction $\alpha$ always ends successfully unless consensus is reached. So, the protocol leads to a series of actions which either continues forever including movements forever or finishing when consensus is reached. We group actions to form a series of updating maps $f_t$. We group by the following rule: Starting with the first action we collect actions in the same group until an agent is found which moves. The next updating map $f_2$ is formed analog starting with the next action, and so on. Thus we have a series of update maps.

It is simple to see that the series of updating maps $f_1,f_2,\dots$ fulfills the conditions of Theorem \ref{thm:familyavmap} with $y$ the identity. Every $f_t$ is an averaging map because by definition the movement of agents goes into the convex hull or along its border and stops before the convex hull is left. It is equiproper, because for each $x$ there are only as many possible updating maps as their convex hull has extreme points. Thus, there is a $\delta(x)>0$ by taking the minimum over this finite set of possible updating maps.
Every $f_t$ is continuous in $x$ when we regard all agents which have the same position as one agent. Equicontinuity at $x$ again follows from finiteness of the possible updating maps. 

Thus, the protocol in Example \ref{ex:watergun} leads to convergence to consensus. This can not be shown by applying Moreaus's Theorem 2 because the movements cannot be easily encoded in terms of communication topologies. One could try to specify it in terms of communication topologies by stating that the moving agent has agents at the detected directions in $A$ as its set of neighbors. But even then the conditions of Moreau's Assumption 1 (especialy number 3) need not be fulfilled and a node connected to all other nodes across time intervals of length $T$ need not exist as necessary for Moreau's Theorem 2. 

\begin{ex}[Nonlinear proper averaging map] Let
  \begin{gather*}
    f_1(x) = x_1, \,  f_2(x) = a(l)x_1 +(1-a(l))x_2,\, f_3(x) = \tfrac{1}{5}x_2 + \tfrac{4}{5}x_3
  \end{gather*}
  where $l = \text{dist}(x_3,\text{line passing through $x_1$ and
    $x_2$})$ 
  and $a$ is continuous and decreasing from $\tfrac{1}{2}$ to 0 in $[0,1]$ and zero otherwise.
In this example agent 3 moves towards agent 2 while agent 2 moves towards agent 1 only if agent 3 is
close to a stripe around the line through agent 2 and agent 1.
\end{ex} 
Examples of this kind can be formulated in terms of communication topologies as Moreau's Theorem 2 needs them, but the existance of a uniform bound for the length of intercommunication intervals $T$ is not easily at hand. 

\begin{ex}[Non-arithmetic means]\label{ex:nonarithmetic}
  We define $g_1,g_2,g_3,g_4: (\R^d)^3 \to \R^d$ by $g_1(x) := \max
  \{x^1,x^2,x^3\}$, $g_2(x) := \frac{1}{3}(x^1+x^2+x^3)$, $g_3(x) :=
  \sqrt[3]{x^1x^2x^3}$ and $g_4(x) := \min\{x^1,x^2,x^3\}$ with all
  computations componentwise. Further on let 
  $f^{\sigma_1\sigma_2\sigma_3}: (\R^d)^3\to(\R^d)^3$ with
  \[
  f^{\sigma_1\sigma_2\sigma_3}:=(g_{\sigma_1},g_{\sigma_2},g_{\sigma_3}).
  \]
\end{ex}
  It is easy to verify, that the family of all
  $f^{\sigma_1\sigma_2\sigma_3}$ where 1 and 4 are not both in
  $(\sigma_1,\sigma_2,\sigma_3)\in\{1,2,3,4\}^3$ is an equicontinuous
  set of $y$-averaging maps, when the $y$-convex hull is the interval
  $[\min_{i\in\n}x^i,\max_{i\in\n}x^i]$. Equiproper is implied by
  finiteness. Thus convergence to consensus is ensured by Theorem
  \ref{thm:familyavmap}. Moreau's theorem is not applicable because
  $f^{\sigma_1\sigma_2\sigma_3}$ is not a convex hull averaging map if
  some $\sigma_i$ is 1, 3 or 4 (since the componentwise $\min$ or
  $\max$ and the geometric mean are in general not contained in the
  convex hull).

Krause~\cite{Krause2009Compromiseconsensusand} shows another example where Moreau's theorem
does not imply convergence: Assume three agents in two dimensional space.
In each iteration every agent takes the
mean value of the two other agents.  Hence, no agent moves into the
relative interior of the convex hull but these maps are still proper
averaging maps and Theorem~\ref{thm:familyavmap} applies.

\appendices
\section{Proofs}
\label{app:proof-theor-familyavmap}

\begin{IEEEproof}[Proof of Lemma~\ref{lem:limit_proper_avmap}]
    First we show that $g$ is an averaging map. Take $x\in S^n$ and let
  $\eps>0$. Due to the pointwise convergence of $(f_{t})_i$ to
  $g_i$ and uniform continuity of $y$ there is $t_0$ such that for all
  $t>t_0$ it holds $\|y^i(f_{t}(x)) - y^i(g(x))\| < \eps$. Due
  to $y^i(f_{t}(x)) \in \yconv(x)$ it follows that the maximal
  distance of $y^i(g(x))$ to $\yconv(x)$ is less than $\eps$, and thus
  $y^i(g(x))\in \yconv(x)$ because $\yconv(x)$ is closed.

  We show that $g$ is proper. To this end, let $x\in S^n$ be not a
  consensus. We have to show that there is $z^\ast \in \yconv(x)$ but
  $z^\ast\notin\yconv(g(x))$. (Note that $z^\ast\in S$, while $x\in
  S^n$ and $y(x)\in S^m$.)  We know that there is for each $t\in\N$ an
  $z(t) \in \yconv(x)$ with $z(t) \notin
  \yconv(f_{t}(x))$. According to the equiproper property it can
  be chosen such that the distance of $z(t)$ to
  $\yconv(f_{t}(x))$ is bigger than $\frac{\delta(x)}{2}>0$ for
  all $t\in\N$. Further on, we know that the set difference
  $\yconv(f_{t}(x))\backslash\yconv(x)$ is non empty and
  bounded, thus there is a subsequence $t_{s}$ such
  that $z(t_s)$ converges to a $z^\ast\in\yconv(x)$. Because of the
  construction it also holds $z^\ast\notin\yconv(g(x))$.
\end{IEEEproof}

\begin{IEEEproof}[Proof of Theorem~\ref{thm:familyavmap}]
  The idea of the proof is the following: We define
  $C(t) := \yconv (x(t))$
  which is convex and
  compact. It holds $C(t+1) \subset C(t)$ because of the averaging
  property and $C:=\bigcap_{t=0}^\infty C(t) \neq \emptyset$ because of
  compactness.
  In the following we will show that $C$ is a singleton, and that
  for all $i\in\n$ the sequences $x^i(t)$ converge to it.
  This will be done in three main steps,
  but first we note that because of compactness of $C(0)^n$ there is a
  subsequence $t_s$ and $c := (c^1,\dots,c^n) \in C(0)^n$ such that
  $\lim_{s\to\infty} x(t_s) = c$.

  \begin{enumerate}
  \item  We show that $C = \yconv(c)$.
    To accept ''$\supset$'' see that for all $t_s\geq t$ there is
    $x^i(t_s)\in C(t)$ and thus $c^i \in C(t)$. This implies $c^i \in C$
    because all the $C(t)$ are closed.

    To show ''$\subset$'' let $x\in C$ and $\eps>0$.
    Because of uniform continuity of $y$ there is $\eta>0$ such that for
    every $x'\in S$ with $\norm{c-x'}<\eta$ it holds
    $\norm{y(c)-y(x')}<\eps$. Further on, there is $s_0$ such that for
    all $s\geq s_0$ it holds $\norm{x(t_s)-c}<\eta$. This implies for
    every $i\in\m$ that $\norm{y^i(x(t_s))-y^i(c)}<\eps$.

    Obviously, $x\in C(t_{s_0})$. Thus, there exist convex coefficients
    $a_1,\dots,a_m\in\R_{\geq 0}^d$ such that $x = \sum_{i=1}^m a_i
    y^i(x(t_{s_0}))$.  Now we can conclude
    \begin{align*}
      \|x-\sum_{i=1}^m a_i y^i(c)\| &= \|\sum_{i=1}^m a_i(y^i(x(t_{s_0}))
      - y^i(c))\| \\ &\leq \sum_{i=1}^m \|y^i(x(t_{s_0})) - y^i(c)\| = m\eps.
    \end{align*}
    It follows that $x \in \yconv(c)$ because
    $\yconv(c)$ is closed.

  \item The next step is to show that $c$ is a consensus, i.e.  $c^1 =
    \dots = c^n$.  The family $F$ is uniformly equicontinuous and for
    all $x\in X$ it holds that $\{f(x) \,|\, f\in F\}$ is bounded (and
    thus relatively compact) because all the $f$ are averaging
    maps. So, due to the theorem of Arzel\`a-Ascoli, $F$ is relatively
    compact.  Thus, there is a subsequence $t_{s_r}$ such that
    $f_{t_{s_r}}$ converges uniformly to a continuous limit function
    $g$ for $r\to\infty$.  Due to Lemma~\ref{lem:limit_proper_avmap}
    we also know that $g$ is a proper averaging map. 
    In two substeps we show that
    $\yconv(g(c)) = \yconv(c)$ which implies that $c$ is a consensus:
    \begin{enumerate}
    \item We show that for all $i\in\n$ it
      holds $\lim_{r\to\infty} f_{t_{s_r}}(x_{t_{s_r}}) = g(c)$. We know
      that $f_{t_{s_r}} \to g$ uniformly and that $x^i(t_{s_r}) \to c$.
      Now we estimate
      \begin{align}
        \|f_{t_{s_r}}(x(t_{s_r})) - g(c)\| \leq&
         \|f_{t_{s_r}}(x(t_{s_r})) -
        f_{t_{s_r}}(c)\| \nonumber \\
        &+ \|f_{t_{s_r}}(c) - g(c)\|\nonumber
      \end{align}
      Both terms on the right hand side can be smaller than
      $\frac{\eps}{2}$ for any $\eps$ for large enough $r$ because of the
      continuity of $f_{t_{s_r}}$ and the uniform convergence
      $f_{t_{s_r}} \to g$.
    
    \item We show
      $\yconv(g(c)) = \yconv(c)$.
      ''$\subset$'' holds because $g$ is an averaging map (see
      2a). To show ''$\supset$'' let
      $x\in\yconv (c)$. Thus, for all $r$ it holds $x\in C
      \subset C(t_{s_r}+1)$ and thus there exist convex coefficients with
      convex combination $x = \sum_{i=1}^{m}a_i(r) y^i(x(t_{s_r}+1))$.
      Now, $(a_1(r), \dots, a_m(r))_{r\in\N}$ is a sequence in the
      compact set of convex coefficients and thus there is a subsequence
      $r_q$ such that $\lim_{q\to\infty}(a_1(r_q), \dots, a_m(r_q)) =
      ({a_1}^\ast \dots {a_m}^\ast)$.  Now due to
      2c and continuity of $y$ it holds,
      \begin{align}
        x&=\sum_{i=1}^{m}\lim_{q\to\infty}a_i(r_q)\lim_{q\to\infty}y^i(x(t_{s_{r_q}}+1)) \nonumber \\
        &= \sum_{i=1}^{m}{a_i}^\ast y^i(g(c)).\nonumber
      \end{align}
      Thus, $x \in \yconv(g(c))$.

    \end{enumerate}
    This implies that $c$ is a consensus, because $g$ is a proper
    averaging map.

  \item Finally, we show that for each
    $i\in\n$ the sequence $(x^i(t))_{t\in\N}$ (and not only subsequences)
    converges to $\gamma := c^1 = \dots = c^n$ for $t\to\infty$.  We know
    that for $\eps > 0$ there is a $r_0$ such that for each $i\in\n$ it
    holds $\|y^i(x(t_{s_{r_0}})) - \gamma\|<\eps$.  Further on, for
    $t\geq t_{s_{r_0}}$ it holds $x(t) \in C(t) \subset C(t_{s_{r_0}})$.
    Thus, for each $i\in\n$ there are convex combinations $x^i(t) =
    \sum_{j=1}^m a^j y^j(x(t_{s_{r_0}}))$.  Now, we conclude for all $t >
    t_{s_{r_0}}$
    \begin{align}
      \|x^i(t) - \gamma\| &= \|\sum_{j=1}^n a^j y^j(x(t_{s_{r_0}})) -
      \gamma)\| \nonumber \\
      &\leq \sum_{j=1}^m \|x^j(t_{s_{r_0}}) -
      \gamma\| = m\eps.\nonumber
    \end{align}
  \end{enumerate}
  This proves the theorem.
\end{IEEEproof}

\begin{IEEEproof}[Proof of Lemma~\ref{lem:cont_dep_initial_value}]
    Let $\eps>0$ and consider two initial values $x(0),\tilde x(0)\in
  S^n$ with corresponding limits $\gamma$, $\tilde\gamma$
  respectively. We have to show that there exists $\delta>0$ such that
  $d(x(0),\tilde x(0))\leq\delta$ implies
  $d(\gamma,\tilde\gamma)\leq\eps$.
  
  We note that for every $t$ it holds that
  \begin{align*}
    d(\gamma,\tilde\gamma) 
    & \leq d(\gamma,x(t)) + d(x(t),\tilde x(t)) + d(\tilde\gamma,\tilde x(t)).
  \end{align*}
  We choose $t_0$ large enough, that
  \begin{align*}
    d(\gamma,x(t_0)) & \leq \frac{\eps}{3} &
    d(\tilde\gamma,\tilde x(t)) & \leq \frac{\eps}{3}.
  \end{align*}
  Since $\{f_t\}$ is an equicontinuous family there exists $\eta>0$
  such that for every $t\in\N$ it holds that
  \[
  d(x(t),\tilde x(t)) \leq \eta\implies d(f_t(x(t)),f_t(\tilde
  x(t)))\leq\eps.
  \]
  Since $x(t)$ and $\tilde x(t)$ solve $x(t+1) = f_t(x(t))$ we have
  recursively that for every  $t_0$ there exists $\delta>0$ such
  that
  \[
  d(x(0),\tilde x(0)) \leq \delta\implies d(x(t_0),\tilde
  x(t_0))\leq\frac{\eps}{3}
  \]
  which implies the claim.
\end{IEEEproof}

\begin{IEEEproof}[Proof of Theorem~\ref{thm:family_yconvphi_avmap}]
  We define $f_t = \phi\circ g_t\circ \phi^{-1}:S^n\to S^n$.  We show that $\{f_t\,|\,
  t\in \N\}$
  is a family of equicontinuous, equiproper $y$-averaging maps on $S^n$.
  Equicontinuity and the fact that the $f_t$'s are $y$-averaging maps are clear.
  To see equiproperness of $f_t$ we note first that from equiproperness of $g_t$ it follows
  \begin{align}
    d_H\bigl(\phi^{-1}(\conv_{i\in\m}y^i(\phi(g_t(x)))),\phi^{-1}(\conv_{i\in\m}y^i(\phi(x)))\bigr) & \geq \delta(x)\nonumber\\\allowdisplaybreaks
    \implies 
    d_H(\phi^{-1}(\yconv(f_t(\xi))),\phi^{-1}(\yconv(\xi))) & \geq \delta(\phi^{-1}(\xi))\label{eq:est_equiprop_ft}\nonumber
  \end{align}
  while the second line holds for all $\xi = \phi(x)\in S^n$ and $t\geq 0$.
  Due to (\ref{eq:d_H_inclusion}) we can express the Hausdorff distance as
  \begin{align*}
    d_H&(\phi^{-1}(\yconv(f_t(\xi))),\phi^{-1}(\yconv(\xi)))  \\\allowdisplaybreaks
    &= \max_{z\in \phi^{-1}(\yconv(f_t(\xi)))}\ \min_{w\in\phi^{-1}(\yconv(\xi))} \norm{z-w}\\\allowdisplaybreaks
    &= \max_{\phi(z)\in \yconv(f_t(\xi))}\ \min_{\phi(w)\in\yconv(\xi)} \norm{z-w}.
  \end{align*}
  With this preparation we show equiproperness of the $f_t$'s:
  \begin{align*}
    d_H&\bigl(\yconv(f_t(\xi)),\yconv(\xi)\bigr) \\
    &= \max_{\zeta\in \yconv(f_t(\xi))}\ \min_{\omega\in\yconv(\xi)} \norm{\zeta-\omega}\\
    &= \max_{\phi(z)\in \yconv(f_t(\xi))}\ \min_{\phi(w)\in\yconv(\xi)} \norm{\phi(z)-\phi(w)} \\
    &\geq L\max_{\phi(z)\in \yconv(f_t(\xi))}\ \min_{\phi(w)\in\yconv(\xi)} \norm{z-w}\\
    &\geq L\delta(\phi^{-1}(\xi))
  \end{align*}
  where $L$ is the Lipschitz constant of $\phi^{-1}$.
  Now for
  $\xi(t) = \phi(x(t))$ it follows
  $%\[
  \xi(t+1) = \phi(g_t(x(t)) = \phi(g_t(\phi^{-1}(\xi(t)))) = f_t(\xi(t)). 
  $ %\]
  By virtue of Theorem~\ref{thm:familyavmap}, $\xi(t)\to c$ where $c\in S^n$ is
  a consensus and hence, $x(t) \to \phi^{-1}(c)\in T^n$ which is also a consensus.
\end{IEEEproof}

\ifCLASSOPTIONcaptionsoff \newpage \fi

\bibliographystyle{IEEEtran}
% \bibliography{refs}

\end{document}